\documentclass[12pt]{amsart}
\usepackage[latin1]{inputenc}
\usepackage{indentfirst}
\usepackage{amsfonts}
\usepackage{graphicx}
\usepackage{amsmath}
\usepackage{amssymb}
\usepackage{latexsym}
\usepackage{amscd}
\usepackage[all]{xy}
\usepackage{xcolor}

\newcommand{\F}{\mathbb {F}}

\newtheorem{theorem}{Theorem}[section]
\newtheorem{definition}[theorem]{Definition}
\newtheorem{lemma}[theorem]{Lemma}
\newtheorem{corollary}[theorem]{Corollary}
\newtheorem{proposition}[theorem]{Proposition}

\begin{document}

\title[Number of $k$-normal elements]{Number of $k$-normal elements over a
finite field}

\author{Josimar J.R. Aguirre and Victor G.L. Neumann}

\maketitle

\vspace{1ex}
\small{Faculdade de Matem\'{a}tica, Universidade Federal de Uberl\^{a}ndia, 
	Av. 
	J. N. 
	\'{A}vila 2121, 38.408-902 Uberl\^{a}ndia -MG, Brazil}

\vspace{8ex}
\noindent
\textbf{Keywords:} Finite fields, normal elements, $k$-normal elements.\\
\noindent
\textbf{MSC:} 12E20, 11T30.

\begin{abstract}

An element $\alpha \in \mathbb{F}_{q^n}$ is a \textit{normal element} over $\mathbb{F}_q$ if the conjugates $\alpha^{q^i}$, $0 \leq i \leq n-1$, are linearly independent over $\mathbb{F}_q$. Hence a \textit{normal basis} for $\mathbb{F}_{q^n}$ over $\mathbb{F}_q$ is of the form $\{\alpha,\alpha^q, \ldots, \alpha^{q^{n-1}}\}$, where $\alpha \in \mathbb{F}_{q^n}$ is normal over $\mathbb{F}_q$. In 2013, Huczynska, Mullen, Panario and Thomson introduce the concept of k-normal elements, as a generalization of the notion of normal elements. In the last few years, several results have been known about these numbers. In this paper, we give an explicit combinatorial formula for the number of $k$-normal elements in the general case, answering an open problem proposed by Huczynska et al. (2013).

\end{abstract}

\maketitle

\section{Introduction}
Let $\mathbb{F}_{q^n}$ be a finite field with $q^n$ elements, where $q$ 
is a prime power and $n$ is a positive integer.  An element $\alpha \in \mathbb{F}_{q^n}$ is called a \textit{normal element} if the set of conjugates $B := \{ \alpha, \alpha^q, \ldots, \alpha^{q^{n-1}}\}$ forms a basis of $\mathbb{F}_{q^n}$ over $\mathbb{F}_q$. In this case B is called a \textit{normal basis} of $\mathbb{F}_{q^n}$ over $\mathbb{F}_q$.

Normal bases are widely used in applications such as cryptography and signal processing due
to the efficiency of exponentiation. In particular, the $q$-th power of the field elements represented using a
normal basis are given by a cyclic shift (see \cite{jung}, \cite{panario} for further details).

An element $\alpha \in \mathbb{F}_{q^n}$ is normal if and only if the polynomials $g_{\alpha}(x)=\alpha x^{n-1}+ \alpha^q x^{n-2} + \ldots + \alpha^{q^{n-2}} x + \alpha^{q^{n-1}}$ and $x^n-1$ are relatively prime over $\mathbb{F}_{q^n}$
\cite[Theorem 2.39]{Nied}. With this as motivation, Huczynska et al. \cite{knormal} introduced the concept of $k$-normal elements, as an extension of the usual definition of normal elements:

\begin{definition}
Let $\alpha \in \mathbb{F}_{q^n}$ and let $g_{\alpha}(x) = \sum_{i=0}^{n-1} \alpha^{q^i}x^{n-1-i} \in \mathbb{F}_{q^n}[x]$. If $\gcd(x^n-1,g_{\alpha}(x))$ over $\mathbb{F}_{q^n}$ has degree $k$ (where $0 \leq k \leq n-1$), then $\alpha$ is a $k$-normal element of $\mathbb{F}_{q^n}$ over $\mathbb{F}_q$.
\end{definition}


From the above definition, elements which are normal in the usual sense are $0$-normal 
and from the Normal Basis Theorem, we know that they always exist \cite[Theorem 2.35]{Nied}.
In \cite[Lemma 3.1]{lucas} the author gives a way to construct $k$-normal elements from a given normal element
and in \cite[Lemma 3.1]{AN} the authors show a numerical condition of the existence of $2$-normal elements.

In \cite{knormal}, the authors established a formula for the number of k-normal elements \cite[Theorem 3.5]{knormal}, using a result due to Ore \cite{ore}. This formula depends on the factorization of $x^n-1$ into irreducibles over $\mathbb{F}_q$, but since the formula obtained is not ``to easy'' to handle numerically, Huczynska et al. proposed the following problem (see \cite[Problem 6.3]{knormal}): For which values of $q, n$ and $k$ can ``nice" explicit formulae (in $q$ and $n$) be obtained for the number of $k$-normal elements of $\mathbb{F}_{q^n}$ over $\mathbb{F}_q$?

On this line, in \cite{STU} the authors obtained some explicit formulas for certain particular cases of $q$ and $n$, the results depend on the explicit factorization of cyclotomic polynomials and the solutions of some linear Diophantine
equations.

In this paper we obtain an explicit formula for all extensions $\mathbb{F}_{q^n}$ over $\mathbb{F}_q$ as follows: In Section 2, we provide background material that is used  throughout the paper. In Section 3, we obtain the explicit formulas for the number of $k$-normal elements for $\mathbb{F}_{q^n}$ over $\mathbb{F}_q$. Finally, in section 4, using the formulas in the previous section,  we obtain formulas for particular cases with $k$ small and we show numerical results for some cases using the software SageMath \cite{SAGE}.

\section{Preliminaries}

In this section, we present some definitions and results that will be useful in this paper.
We start with the following definition.

\begin{definition}
Let $f(x)$ be a monic polynomial with coefficients in $\mathbb{F}_{q}$. The Euler Totient Function for polynomials over $\mathbb{F}_q$ is given by
		$$
		\Phi_q(f)= \left| \left( \dfrac{\mathbb{F}_q[x]}{\langle f \rangle} \right)^{*} \right|,
		$$
		where $\langle f \rangle$ is the ideal generated by $f(x)$ in $\mathbb{F}_q[x]$.
\end{definition}

From an equivalent definition of the notion of $k$-normal elements  (see \cite[Theorem 3.2]{knormal}) and
using the Euler Totient Function for polynomials above,
we get the formula for the number of $k$-normal elements.

\begin{theorem}\label{counting-knormal}
\cite[Theorem 3.5]{knormal}
The number of $k$-normal elements of $\mathbb{F}_{q^n}$ over $\mathbb{F}_q$ is given by
	$$
	\sum_{\substack{h| x^n-1 \\ \deg(h)=n-k}} \Phi_q(h),
	$$
	where the divisors are monic and the polynomial division is over $\mathbb{F}_q$.
\end{theorem}

In the last formula, we see that the knowledge of the irreducible factors of $x^n-1$ is very important. In 
that sense, the following results will help us to obtain ``nice'' formulas for the number of $k$-normal elements.

\begin{lemma}\label{divisoresxn-1} Let $q$ be a prime power and 
let $n$ be a positive integer.
Let $v_d$ be the number of distinct irreducible monic factors of $x^n-1$ with degree $d$ over $\F_q$. We have
\begin{equation}\label{factors}
v_d= \dfrac{1}{d} \sum_{r \mid d} t_r 
\mu \left( \frac{d}{r} \right),
\end{equation}
where $t_r:=\gcd(q^r-1,n)$.
\end{lemma}
\begin{proof}
Let $\alpha$ be a primitive element in $\F_{q^d}$. 
The number of elements $\alpha^s$, with $1 \leq s \leq q^d-1$, 
with $(\alpha^s)^n=1$ is 
$t_d$. 
Observe that for each irreducible polynomial 
of degree $d$ 
defined over $\F_{q}$ which divides $x^n-1$, there are $d$ elements $\alpha$ in $\F_{q^d}$ such that $\alpha \notin \F_{q^r}$ for any $r<d$, with $\alpha^n=1$. So, from the 
definition of $v_d$, the number of elements $\alpha$ in $\F_{q^d}$
which are not in $\F_{q^r}$ for $r<d$, with $\alpha^n=1$  is
$d 
v_d$, so
$t_d = \sum_{r \mid d} r v_r$.
From the Möbius inversion formula (see \cite[Theorem 2.9]{Apostol}) we get
the desired result.
\end{proof}

Let $p$ be the characteristic of $\mathbb{F}_q$ and 
let $n=p^s n_0$ such that $\gcd (n_0 , p ) = 1$. 
Observe that $x^n-1 = (x^{n_0}-1)^{p^s}$ and the number of distinct irreducible monic factors of $x^n-1$ and $x^{n_0}-1$ with degree $d$ over $\F_q$ are equal.

\begin{corollary}
Let $p$ be the characteristic of $\mathbb{F}_q$, 
let $n=p^s  n_0$ such that $\gcd (n_0 , p ) = 1$ and
let $d$ be the least positive integer such that $n_0 \mid q^d -1$.
The number of different irreducible factors of $x^n-1$ over $\mathbb{F}_q$ is
$$
\frac{1}{d}\sum_{r \mid d} \gcd(q^r-1,n) \cdot \phi\left( \frac{d}{r} \right) ,
$$
where $\phi$ is the Euler totient function.
\end{corollary}
\begin{proof}
Let $\alpha$ be a root of $x^n-1=0$. 
Since $x^n-1=(x^{n_0}-1)^{p^s}$ we get 
$\alpha \in \mathbb{F}_{q^d}$. This implies that if
$r$ is the least positive integer such that $\alpha \in \mathbb{F}_{q^r}$, then $r \mid d$.
Thus, if $r \nmid d$ then $v_r=0$. Let $\omega_q(x^n-1)$ be the number of different irreducible
factors of $x^n-1$ over $\mathbb{F}_q$. From considerations above and Lemma \ref{divisoresxn-1} we have
$$
\omega_q(x^n-1) =
\sum_{r \mid d} v_r = 
\sum_{r\mid d}
\left(
\frac{1}{r}\sum_{u \mid r} t_u \cdot \mu \left( \frac{r}{u} \right)
\right).
$$
Writing $r=uv$ we get
$$
\omega_q(x^n-1) =
\sum_{u \mid d} \sum_{v \mid \frac{d}{u}}
\frac{1}{uv} \cdot t_u \cdot \mu ( v )
=
\sum_{u \mid d} \frac{t_u}{u} \sum_{v \mid \frac{d}{u}}
\frac{\mu(v)}{v} .
$$
From \cite[Theorem 2.3]{Apostol} we get that 
$\sum_{v \mid \frac{d}{u}} \frac{\mu(v)}{v} = \frac{\phi(\frac{d}{u})}{\frac{d}{u}}$, so
$$
\omega_q(x^n-1) =
\sum_{u \mid d}  \frac{t_u \cdot \phi(\frac{d}{u})}{d} =
\frac{1}{d}\sum_{u \mid d} t_u \cdot \phi\left( \frac{d}{u} \right) .
$$

\end{proof}

\section{The number of $k$-normal elements}

We will start this section by finding a formula for the number of normal elements over finite fields. It is worth mentioning that similar formulas already exist in the literature (see \cite[Corollary 5.2.8]{panario}), but we will show them here to observe the ideas that will be used in the general case.

\begin{theorem}\label{normals}
Let $p$ be the characteristic of $\mathbb{F}_q$ and 
let $n=p^s  n_0$ such that $\gcd (n_0 , p ) = 1$.
The number of normal elements of
$\mathbb{F}_{q^n}$ over $\mathbb{F}_q$
is
\begin{equation}\label{numbernormals}
N_0 := q^{n- n_0 }\prod_{r\mid d} (q^r - 1)^{v_ r} ,
\end{equation}
where 
$d$ is the least positive integer such that $q^d \equiv 1 \pmod {n_0}$. 
\end{theorem}
\begin{proof}
From Theorem \ref{counting-knormal}, the number of normal elements of
$\mathbb{F}_{q^n}$ over $\mathbb{F}_q$ is $\Phi_q (x^n -1)$.

Since $q^d \equiv 1 \pmod {n_0}$, 
all the roots
of $x^{n_0}-1$ are in $\mathbb{F}_{q^{d}}$. Thus, if $\xi$ is a root of $x^{n_0}-1$ and $r$ is the smallest positive integer such that
$\xi \in \mathbb{F}_{q^r}$, then $r \mid d$. 
This means that if $v_r \neq 0$ then $r \mid d$.
Let $g_{r,1},\ldots , g_{r,v_r}$ be the $v_r$ irreducible monic polynomials of degree $r$ which divide $x^{n_0} -1$. Thus
$$
x^n -1 =
\left( \prod_{r \mid d} \prod_{i=1}^{v_r} g_{r,i}
\right)^{p^s}.
$$
Since $\Phi_q((g_{r,i})^{p^s})=(q^r)^{p^s} - (q^r)^{p^s - 1} = q^{r(p^s - 1)} (q^r - 1)$ then
$$
\Phi_q (x^n -1)  = \prod_{r\mid d} q^{r v_r (p^s-1)} (q^r - 1)^{v_ r} .
$$
From Lemma \ref{divisoresxn-1}, we have
$$
\sum_{r\mid d} r v_r = 
\sum_{r\mid d} \sum_{u \mid r} t_u \cdot \mu \left( \frac{r}{u} \right)
= \sum_{u \mid d} \sum _{k \mid \frac{d}{u}} t_u \cdot \mu (k)
= \sum_{u \mid d} t_u \sum_{k \mid \frac{d}{u}}  \mu (k)=t_d=n_0.
$$

%
%

Substituting that into the formula of $\Phi_q(x^n-1)$ 
and taking into account that $n_0(p^s-1)= n - n_0$, we get
$$
\Phi_q (x^n -1)  = q^{n- n_0 }\prod_{r\mid d} (q^r - 1)^{v_ r} .
$$
\end{proof}
With equation \eqref{numbernormals}  we may calculate easily the number of normal elements of 
$\mathbb{F}_{q^n}$ over $\mathbb{F}_q$ with the use of SageMath (see \cite{SAGE}).
To find the value of $d$ remember that $d \mid \phi(n_0)$.
We illustrate this formula in Tables \ref{table1}, \ref{table2} and \ref{table3}.

\begin{table}[h]
	\centering
	\begin{tabular}{c|cccccccccc}
		$n$        & $1$ & $2$ & $3$  &$4$ & $5$    & $6$   & $7$    & $8$ & $9$ & $10$\\
		$N_0$  &  $1$ & $2$ & $3$  &$8$ & $15$ & $24$ & $49$ & $128$ & $189$ & $480$\\
		\hline		\hline
$n$        &  $11$  &$12$             & $13$       & $14$   & $15$           & $16$ & $17$                  & $18$ & $19$                     & $20$ \\
$N_0$  &  $1023$ & $1536$ & $4095$  &$6272$ & $10125$ & $ 32768$ & $65025$ & $96768$ & $262143$ & $491520$ \\
	\end{tabular}\vspace*{0.5cm}
	\caption{The number of normal elements of $\mathbb{F}_{q^n}$
over $\mathbb{F}_q$ with $q=2$.}
	\label{table1}
\end{table}

\begin{table}[h]
	\centering
	\begin{tabular}{c|cccccccc}
		$n$        & $1$ & $2$ & $3$  &$4$ & $5$    & $6$   & $7$    & $8$ \\
		$N_0$  &  $2$ & $4$ & $18$  &$32$ & $160$ & $324$ & $1456$ & $2048$ \\
		\hline		\hline
		$n$        &  $9$& $10$                   & $11$           &$12$             & $13$            & $14$              & $15$           & $16$ \\
		$N_0$  &  $13122$ & $25600$ & $117128$  &$209952$ & $913952$ & $2119936$ & $9447840$ & $13107200$ \\
	\end{tabular}\vspace*{0.5cm}
	\caption{The number of normal elements of $\mathbb{F}_{q^n}$
over $\mathbb{F}_q$ with $q=3$.}
	\label{table2}
\end{table}

\begin{table}[h]
	\centering
	\begin{tabular}{c|ccccccc}
$n$        & $1$ & $2$ & $3$        &$4$      & $5$      & $6$       & $7$         \\
$N_0$  &  $3$ & $12$ & $27$  &$192$ & $675$ & $1728$ & $11907$  \\
\hline		\hline
$n$         & $8$ &  $9$             & $10$                   & $11$           &$12$             & $13$            & $14$       \\
$N_0$  & $49152$ &  $107163$ & $691200$ & $3139587$  &$7077888$ & $50307075$ & $195084288$ \\
\end{tabular}\vspace*{0.5cm}
\caption{The number of normal elements of $\mathbb{F}_{q^n}$
over $\mathbb{F}_q$ with $q=4$.}
\label{table3}
\end{table}

%
%

Let us now look at the general case of $k$-normal elements, with a generalization of the previous idea.

\begin{theorem}\label{numberknormals}
Let $p$ be the characteristic of $\mathbb{F}_q$ and 
let $n=p^s  n_0$ such that $\gcd (n_0 , p ) = 1$.
Let $k$ be a non-negative integer such that $k <n$ and 
$d$ be the least positive integer such that $q^d \equiv 1 \pmod {n_0}$.
The number of $k$-normal elements
of $\mathbb{F}_{q^n}$ over $\mathbb{F}_q$ is
\begin{equation}\label{numknormales}
N_k := \sum_{(\alpha_{r,i}) \in \mathcal{A}_d}  \prod_{r \mid d}
\prod_{\substack{i=1 \\ \alpha_{r,i}>0}}^{v_r} q^{r (\alpha_{r,i}-1)} (q^r-1).
\end{equation}
where $\mathcal{A}_d$ is the set of tuples $(\alpha_{r,i})$, for $r \mid d$ and $1 \leq i \leq v_r$,
such that $0 \leq \alpha_{r,i} \leq p^s$ and
$$
\sum_{r\mid d} r  \sum_{i=1}^{v_r} \alpha_{r,i} = n-k .
$$
In particular, if $n$ is prime to $q$ (i.e. $n=n_0$) then
$$
N_k = \sum_{(a_r)\in A_d} \prod_{r\mid d} \binom{v_r}{a_r} (q^r - 1)^{a_r},
$$
where $A_d$ is the set of tuples $(a_r)_{r \mid  d}$ such that
$\displaystyle \sum_{r\mid d} ra_r = n-k$ and $0 \leq a_r \leq v_r$.
\end{theorem}
\begin{proof}
From Theorem \ref{counting-knormal}, the number of  $k$-normal elements of $\mathbb{F}_{q^n}$ over $\mathbb{F}_q$ is
	$$
	\sum_{\substack{h| x^n-1 \\ \deg(h)=n-k}} \Phi_q(h) .
	$$
From the proof of Theorem \ref{normals}, we have
$$
x^n -1 =
\left( \prod_{r \mid d} \prod_{i=1}^{v_r} g_{r,i}
\right)^{p^s}.
$$
If $h$ is a monic factor of  $x^n-1$ of degree  $n-k$ then
\begin{equation}\label{divn_k}
h = \prod_{r \mid d} \prod_{i =1}^{v_r} g_{r,i}^{\alpha_{r,i}},
\end{equation}
where $(\alpha_{r,i}) \in \mathcal{A}_d$.
For this polynomial we have
$$
\Phi_q(h) = \prod_{r \mid d} \prod_{\substack{i=1 \\ \alpha_{r,i}>0}}^{v_r} q^{r (\alpha_{r,i}-1)} (q^r-1),
$$
so 
$$
N_k = \sum_{(\alpha_{r,i}) \in \mathcal{A}_d}  \prod_{r \mid d}
      \prod_{\substack{i=1 \\ \alpha_{r,i}>0}}^{v_r} q^{r (\alpha_{r,i}-1)} (q^r-1).
$$
In the case where $n=n_0$ (i.e. $s=0$ ) each $\alpha_{r,i}$ is $0$ or $1$.
Let $A_d$ be the set of tuples $(a_r)_{r \mid  d}$ such that
$\displaystyle \sum_{r\mid d} ra_r = n-k$ and $0 \leq a_r \leq v_r$. For each 
$a_r$ there are $\displaystyle \binom{v_r}{a_r}$ 
choices of $(\alpha_{r,1},\ldots , \alpha_{r,v_r})$ with
$\sum_{i=1}^{v_r} \alpha_{i,v_i} = a_r$. For each fixed choice of 
$(a_r)_{r \mid  d}$ and $(\alpha_{r,1},\ldots , \alpha_{r,v_r})$
there is a polynomial $h$ of the form \eqref{divn_k}, and for this polynomial we have
$$
\Phi_q(h) = \prod_{r \mid d} \prod_{\substack{i=1 \\ \alpha_{r,i}>0}}^{v_r} q^{r (\alpha_{r,i}-1)} (q^r-1)
=\prod_{r \mid d} (q^r-1)^{a_r}.
$$
Substituing $\Phi_q(h)$ in the formula of $k$-normal elements (see Theorem\ref{counting-knormal}) we get
$$
N_k = \sum_{(a_r)\in A_d} \prod_{r\mid d} \binom{v_r}{a_r} (q^r - 1)^{a_r}.
$$
\end{proof}

The combinatorial formula obtained to find the number of $k$-normal elements 
allows finding $N_k$ for particular values of $q$ and $n$ in a simpler way, 
without depending on the factoring of cyclotomic polynomials as in Theorem \ref{counting-knormal}. 
For example, in \cite{STU} the authors found the value of $N_k$ for $n=p^n$ and $q \equiv 1 \pmod p$, 
where $p$ is a prime number (see \cite[Theorem 8]{STU}), which can be obtained directly from the theorem above.

In \cite{Tinani}, the authors found a lower bound on the number of $k$-normal elements
in $\mathbb{F}_{q^n}$ when they exist. We retrieve the same result using \eqref{numknormales}.

\begin{proposition}\cite[Theorem 4]{Tinani}
Let $k\in \{0,1,\ldots , n \}$ and let $N_k$ the number of $k$-normal elements in
$\mathbb{F}_{q^n}$. If $N_k>0$, i.e. if $k$-normal elements exists in $\mathbb{F}_{q^n}$,
then
$$
N_k \geq \frac{N_0}{q^k}.
$$
\end{proposition}
\begin{proof}
We have
$$
N_k =  \sum_{(\alpha_{r,i}) \in \mathcal{A}_d}  \prod_{r \mid d}
\prod_{\substack{i=1 \\ \alpha_{r,i}>0}}^{v_r} q^{r (\alpha_{r,i}-1)} (q^r-1) 
\\
 \geq \displaystyle \sum_{(\alpha_{r,i}) \in \mathcal{A}_d} 
 \displaystyle \prod_{r \mid d}
	\prod_{i=1}^{v_r} q^{r (\alpha_{r,i}-1)} (q^r-1).
$$
Since
$$
\sum_{r \mid d} \sum_{i=1}^{v_r} r(\alpha_{r,i}-1) =
n - k - \sum_{r \mid d} rv_r  = n - k - \gcd (q^d-1,n_0)
=n-k-n_0,
$$
and
$\displaystyle
\prod_{r \mid d} (q^r-1)^{v_r} = \frac{N_0}{q^{n-n_0}}$,
we get
$N_k \geq \displaystyle
|\mathcal{A}_d| q^{n-k-n_0} \frac{N_0}{q^{n-n_0}}
\geq \frac{N_0}{q^k}$.
\end{proof}

\section{Number of $k$-normal elements for $k$ small}

In this section we will show the formula of the number of $k$-normal elements in $\mathbb{F}_{q^n}$ over $\mathbb{F}_q$  for $k$ small ($k=1,2,3$), and we will give the exact amount of theses numbers for particular values of $q$,
using SageMath.

\begin{proposition}\label{1normals}
Let $p$ be the characteristic of $\mathbb{F}_q$ and 
let $n=p^s  n_0$ such that $\gcd (n_0 , p ) = 1$.
Let $d$ be the least positive integer such that $q^d \equiv 1 \pmod {n_0}$.
The number of $1$-normal elements
of $\mathbb{F}_{q^n}$ over $\mathbb{F}_q$ is
\begin{equation}\label{number1normal}
N_1 = 
\left\{
\begin{array}{ll}
v_1 q^{n- n_0 - 1}\displaystyle \prod_{r\mid d} (q^r - 1)^{v_ r}  & \text{if } n \neq n_0 ;\\
v_1 (q-1)^{v_1 -1 }\displaystyle \prod_{\substack{r \mid d \\ r \neq 1}} (q^r - 1)^{v_ r}   & \text{if } n = n_0 .
\end{array}
\right.
\end{equation}
\end{proposition}
\begin{proof}
Supppose first that $n\neq n_0$ (so $p^s>1$).
The set $\mathcal{A}_d$, like in Theorem \ref{numberknormals}, has
$v_1=\gcd(q-1,n)$ elements.  
For every $j \in \{1,\ldots , v_1\}$ there is an element $(\alpha_{r,i}) \in \mathcal{A}_d$ such that
$\alpha_{r,i}=p^s$ for $(r,i)\neq (1,j)$ and
$\alpha_{1,j}=p^s-1$. From Theorem \ref{numberknormals}, we get
$$
N_1= \frac{v_1}{q} \prod_{r \mid d}
\prod_{i=1}^{v_r} q^{r (p^s-1)} (q^r-1)
=v_1 q^{n- n_0 - 1}\displaystyle \prod_{r\mid d} (q^r - 1)^{v_ r}.
$$
Supppose now that $n= n_0$.
The set $A_d$ has only  one
element $(a_r)_{r \mid  d} \in A_d$ such that
for $r\neq 1$ we have  $a_r=v_r$ and for $r=1$ we have $a_1=v_1-1$.
From Theorem \ref{numberknormals}, we get
$$
N_1 
=v_1 (q-1)^{v_1 -1 }\displaystyle \prod_{\substack{r \mid d \\ r \neq 1}} (q^r - 1)^{v_ r} .
$$
\end{proof}

Comparing Theorem \ref{normals} with Proposition \ref{1normals}, we get that for $n\neq n_0$
we have $N_1 = \dfrac{v_1}{q} N_0$ and for
$n = n_0$
we have $N_1 = \dfrac{v_1}{q-1} N_0$.

\begin{proposition}\label{2normals}
Let $p$ be the characteristic of $\mathbb{F}_q$ and 
let $n=p^s n_0$ such that $\gcd (n_0 , p ) = 1$.
Let $d$ be the least positive integer such that $q^d \equiv 1 \pmod {n_0}$.
The number of $2$-normal elements
of $\mathbb{F}_{q^n}$ over $\mathbb{F}_q$ is
	\begin{equation}\label{number2normal}
		N_2 = 
		\left\{
		\begin{array}{ll}
\left( v_1 +
\frac{1}{2} v_1(v_1-1) +
v_2
\right)	
q^{n-n_0 - 2}
\displaystyle \prod_{r \mid d} (q^r-1)^{v_r}   & \text{if } n \neq n_0 \text{ and } p^s>2;\\
\left(\frac{q}{q-1} v_1 +
\frac{1}{2} v_1(v_1-1) +
v_2
\right)	
q^{n-n_0 - 2}
\displaystyle \prod_{r \mid d} (q^r-1)^{v_r}   & \text{if } n \neq n_0 \text{ and } p^s=2;\\
\left(
\dfrac{v_1(v_1-1)}{2(q-1)^2} + \dfrac{v_2}{q^2-1}
\right)\displaystyle \prod_{r \mid d} (q^r - 1)^{v_ r}  & \text{if } n = n_0 .
		\end{array}
		\right.
	\end{equation}
\end{proposition}
\begin{proof}
Supppose first that $n\neq n_0$ (so $p^s>1$).
The set $\mathcal{A}_d$ in Theorem \ref{numberknormals} is composed by three kind of elements, 
since
$p^s  \sum_{r\mid d} r  v_r= n$.

a) For every $j \in \{1,\ldots , v_1\}$ there is an element $(\alpha_{r,i}) \in \mathcal{A}_d$ such that
$\alpha_{r,i}=p^s$ for $(r,i)\neq (1,j)$ and 
$\alpha_{1,j}=p^s-2$ (we must consider two cases: $p^s =2$ and $p^s>2$).

b) For every $\{j_1 \neq j_2\} \subset \{1,\ldots , v_1\}$ there is an element $(\alpha_{r,i}) \in \mathcal{A}_d$ such that
$\alpha_{r,i}=p^s$ for $(r,i)\notin \{ (1,j_1), (1,j_2) \}$ and
$\alpha_{1,j_1}=\alpha_{1,j_2}=p^s-1$ (only possible if $v_1>1$).

c) For every $j \in \{1,\ldots , v_2\}$ there is an element $(\alpha_{r,i}) \in \mathcal{A}_d$ such that
$\alpha_{r,i}=p^s$ for $(r,i)\neq (2,j)$ and 
$\alpha_{2,j}=p^s-1$ (only possible if $v_2>0$).

Observe that $v_2>0$ is equivalent to $\gcd(q^2-1,n)>\gcd(q-1,n)$, which is also equivalent
to $\gcd(q+1,n)>2$, or $\gcd(q+1,n)=2$ and
$\gcd(q-1,n) \mid \frac{n}{2}$.

If $p^s>2$,
from Theorem \ref{numberknormals} we get
\begin{eqnarray*}
N_2 & = & 
\left( \frac{v_1}{q^2}  +
\frac{v_1(v_1-1)}{2q^2}  +
	\frac{v_2}{q^2} 
\right)	
q^{n-n_0 }\prod_{r \mid d} (q^r-1)^{v_r} \\
& = & 
\left( v_1 +
\frac{1}{2} v_1(v_1-1) +
v_2
\right)	
q^{n-n_0 - 2}\prod_{r \mid d} (q^r-1)^{v_r} .
\end{eqnarray*}

If $p^s=2$,
from Theorem \ref{numberknormals} we get
\begin{eqnarray*}
	N_2 & = & 
	\left( \frac{v_1}{q(q-1)}  +
	\frac{v_1(v_1-1)}{2q^2}  +
	\frac{v_2}{q^2} 
	\right)	
	q^{n-n_0 }\prod_{r \mid d} (q^r-1)^{v_r} \\
	& = & 
	\left( \frac{q}{q-1}v_1 +
	\frac{1}{2} v_1(v_1-1) +
	v_2
	\right)	
	q^{n-n_0 - 2}\prod_{r \mid d} (q^r-1)^{v_r} .
\end{eqnarray*}

Supppose now that $n= n_0$.
The set $A_d$ has at most two
elements. If $v_1 \geq 2$, there is an element
$(a_r)_{r \mid  d} \in A_d$ such that
for $r\neq 1$ we have  $a_r=v_r$ and for $r=1$ we have $a_1=v_1-2$.
If $v_2 \geq 1$, there is an element
$(a_r)_{r \mid  d} \in A_d$ such that
for $r\neq 2$ we have  $a_r=v_r$ and for $r=2$ we have $a_1=v_2-1$.
From Theorem \ref{numberknormals}, we get
\begin{eqnarray*}
N_2 & = &
\dfrac{v_1(v_1-1)}{2(q-1)^2} \displaystyle \prod_{r \mid d} (q^r - 1)^{v_ r} +
\dfrac{v_2}{q^2-1} \displaystyle \prod_{r \mid d} (q^r - 1)^{v_ r} \\
  & = &
\left(
\dfrac{v_1(v_1-1)}{2(q-1)^2} + \dfrac{v_2}{q^2-1}
\right)\displaystyle \prod_{r \mid d} (q^r - 1)^{v_ r} .
\end{eqnarray*}
Observe that if $v_1=1$ or $v_2=0$ this formula is also true.
\end{proof}

\begin{proposition}\label{3normals}
	Let $p$ be the characteristic of $\mathbb{F}_q$ and 
	let $n=p^s n_0$ such that $\gcd (n_0 , p ) = 1$.
	Let $d$ be the least positive integer such that $q^d \equiv 1 \pmod {n_0}$.
	The number of $3$-normal elements
	of $\mathbb{F}_{q^n}$ over $\mathbb{F}_q$ is
	\begin{equation}\label{number3normal}
		N_3 = 
		\left\{
		\begin{array}{ll}
\left( v_1 + \frac{1}{6} v_1(v_1-1)(v_1+4) + v_1 v_2 + v_3
\right)	
q^{n-n_0 - 3}\displaystyle \prod_{r \mid d} (q^r-1)^{v_r} & \text{if } p^s>3;\\
\left( \frac{q}{q-1} v_1 + \frac{1}{6} v_1(v_1-1)(v_1+4) + v_1 v_2 + v_3 \right)	
q^{n-n_0 - 3}\displaystyle \prod_{r \mid d} (q^r-1)^{v_r}   & \text{if } p^s=3;\\
\left( \left( \frac{q}{q-1} +\frac{v_1 - 2}{6}  \right) v_1(v_1-1) + v_1 v_2 + v_3
\right)	
q^{n-n_0 - 3}\displaystyle \prod_{r \mid d} (q^r-1)^{v_r}   & \text{if } p^s=2;\\
\left(
\frac{v_1(v_1-1)(v_1-2)}{6(q-1)^3} + \frac{v_1v_2}{(q-1)(q^2-1)} +\frac{v_3}{q^3-1}
\right)
\displaystyle \prod_{r \mid d} (q^r - 1)^{v_ r}  & \text{if } n = n_0 .
		\end{array}
		\right.
	\end{equation}
\end{proposition}
\begin{proof}
Supppose first that $n\neq n_0$ (so $p^s>1$).
The set $\mathcal{A}_d$ in Theorem \ref{numberknormals} 
is composed by four kind of elements.
	
a) For every $j \in \{1,\ldots , v_1\}$ there is an element $(\alpha_{r,i}) \in \mathcal{A}_d$ such that
$\alpha_{r,i}=p^s$ for $(r,i)\neq (1,j)$ and 
$\alpha_{1,j}=p^s-3$ (we must consider two cases: $p^s =3$ and $p^s>3$).
	
b) For every $\{j_1 , j_2\} \subset \{1,\ldots , v_1\}$ there is an element $(\alpha_{r,i}) \in \mathcal{A}_d$ such that
$\alpha_{r,i}=p^s$ for $(r,i)\notin \{ (1,j_1), (1,j_2) \}$,
$\alpha_{1,j_1}=p^s-2$ and $\alpha_{1,j_2}=p^s-1$ (there are two cases: $p_s=2$ and $p_s>2$).

c) For every $\{j_1 , j_2, j_3 \} \subset \{1,\ldots , v_1\}$ there is an element $(\alpha_{r,i}) \in \mathcal{A}_d$ such that
$\alpha_{r,i}=p^s$ for $(r,i)\notin \{ (1,j_1), (1,j_2) ,(1,j_3)\}$ and
$\alpha_{1,j_1}=\alpha_{1,j_2}=\alpha_{1,j_3}=p^s-1$ (there are $\frac{v_1(v_1-1)(v_1-2)}{6}$ elements of this kind).

d) For every $j_1 \in \{1,\ldots , v_1\}$ and $j_2 \in \{1,\ldots , v_2\}$ there is an element $(\alpha_{r,i}) \in \mathcal{A}_d$ such that
$\alpha_{r,i}=p^s$ for $(r,i)\notin \{ (1,j_1), (2,j_2) \}$,
$\alpha_{1,j_1}=p^s-1$ and $\alpha_{2,j_2}=p^s-1$.

e) For every $j \in \{1,\ldots , v_3\}$ there is an element $(\alpha_{r,i}) \in \mathcal{A}_d$ such that
$\alpha_{r,i}=p^s$ for $(r,i)\neq (3,j)$ and 
$\alpha_{3,j}=p^s-1$.

If $p^s>3$, from Theorem \ref{numberknormals} we get
\begin{eqnarray*}
N_3 & = & 
\left( \frac{v_1}{q^3}  +
		\frac{v_1(v_1-1)}{q^3}  + \frac{v_1(v_1-1)(v_1-2)}{6q^3}+
		\frac{v_1 v_2}{q^3} + \frac{v_3}{q^3}
		\right)	
		q^{n-n_0 }\prod_{r \mid d} (q^r-1)^{v_r} \\
		& = & 
		\left( v_1 +\frac{1}{6} v_1(v_1-1)(v_1+4)+ v_1 v_2 + v_3
		\right)	
		q^{n-n_0 - 3}\prod_{r \mid d} (q^r-1)^{v_r} .
	\end{eqnarray*}
If $p^s=3$, 
from Theorem \ref{numberknormals} we get
\begin{eqnarray*}
	N_3 & = & 
	\left( \frac{v_1}{q^2(q-1)}  +
	\frac{v_1(v_1-1)}{q^3}  +\frac{v_1(v_1-1)(v_1-2)}{6q^3}+
	\frac{v_1 v_2}{q^3} + \frac{v_3}{q^3}
	\right)\cdot 	\\&&
	q^{n-n_0 }\prod_{r \mid d} (q^r-1)^{v_r} \\
	& = & 
	\left( \frac{q}{q-1} v_1 +  \frac{1}{6} v_1(v_1-1)(v_1+4)+  v_1 v_2 + v_3
	\right)	
	q^{n-n_0 - 3}\prod_{r \mid d} (q^r-1)^{v_r} .
\end{eqnarray*}
If $p^s=2$, 
from Theorem \ref{numberknormals} we get
\begin{eqnarray*}
	N_3 & = & 
	\left( 
	\frac{v_1(v_1-1)}{q^2(q-1)}  + \frac{v_1(v_1-1)(v_1-2)}{6q^3}+
	\frac{v_1 v_2}{q^3} + \frac{v_3}{q^3}
	\right)	
	q^{n-n_0 }\prod_{r \mid d} (q^r-1)^{v_r} \\
	& = & 
	\left( \left( \frac{q}{q-1} +\frac{v_1 - 2}{6}  \right) v_1(v_1-1) + v_1 v_2 + v_3
	\right)	
	q^{n-n_0 - 3}\prod_{r \mid d} (q^r-1)^{v_r} .
\end{eqnarray*}

Supppose now that $n= n_0$.
The set $A_d$ has at most three
elements. If $v_1 \geq 3$, there is an element
$(a_r)_{r \mid  d} \in A_d$ such that
for $r\neq 1$ we have  $a_r=v_r$ and for $r=1$ we have $a_1=v_1-3$.
If $v_2 \geq 1$, there is an element
$(a_r)_{r \mid  d} \in A_d$ such that
for $r>2$ we have  $a_r=v_r$ and for $r=1$ and $r=2$ we have $a_1=v_1-1$ and $a_2=v_2-1$.
If $v_3 \geq 1$, there is an element
$(a_r)_{r \mid  d} \in A_d$ such that
for $r\neq 3$ we have  $a_r=v_r$ and for $r=3$ we have $a_3=v_3-1$.

From Theorem \ref{numberknormals}, we get
$$
N_3 =
\left(
\dfrac{v_1(v_1-1)(v_1-2)}{6(q-1)^3} + \dfrac{v_1v_2}{(q-1)(q^2-1)}
+\dfrac{v_3}{q^3-1}
\right)\displaystyle \prod_{r \mid d} (q^r - 1)^{v_ r} .
$$
\end{proof}


We illustrate the formulas \eqref{numbernormals}, \eqref{number1normal}, \eqref{number2normal} and
\eqref{number3normal} in Tables \ref{table4}, \ref{table5} and \ref{table6} using the software Sagemath. 
\begin{table}[h]
\centering
\begin{tabular}{c|ccccccc}
$n$        & $1$ & $2$ & $3$        &$4$      & $5$       & $6$         & $7$         \\
\hline
$N_0$  &  $24$ & $576$ & $13824$  &$331776$   & $9375000$ & $191102976$ & $5858625024$  \\
$N_1$  &  $1$ & $48$  & $1728$   &$55296$     & $375000$  & $47775744$  & $244109376$  \\
$N_2$  &  $0$ & $1$  & $72$       &$3456$     & $1500$    & $ 4976640$  & $0$  \\
$N_3$  &  $0$ & $0$   & $1$       &$96$       & $600$     & $276480$    & $749952$ 
\end{tabular}\vspace*{0.5cm}
\caption{The number of $k$-normal elements of $\mathbb{F}_{q^n}$
over $\mathbb{F}_q$ with $0 \leq k \leq 3$ and $q=5^2$.}
	\label{table4}
\end{table}

\begin{table}[h]
\centering
\begin{tabular}{c|cccc}
$n$    & $9$             & $10$              &$11$               & $12$                \\
\hline
$N_0$  & $7343167948506$ & $190921648153600$ &$5353168688317736$ & $138991482929321568$  \\
$N_1$  & $271969183278$  & $14686280627200$  &$205891103396836$  & $10295665402171968$   \\
$N_2$  & $10072932714$   & $282428473600$    &$0$                & $762641881642368$      \\
$N_3$  & $373071582$     & $0$               &$0$                & $42912185647968$      
\end{tabular}\vspace*{0.5cm}
\caption{The number of $k$-normal elements of $\mathbb{F}_{q^n}$
over $\mathbb{F}_q$ with $0 \leq k \leq 3$ and $q=3^3$.}
\label{table5}
\end{table}

\begin{table}[h]
\centering
\begin{tabular}{c|ccc}
$n$    & $14$                & $15$                 &$16$                   \\
\hline
$N_0$  & $67521013088256000$ & $437893890380859375$ &$17293822569102704640$ \\
$N_1$  & $4220063318016000$  & $437893890380859375$ &$1080863910568919040$  \\
$N_2$  & $281337554534400$   & $204350482177734375$ &$67553994410557440$      \\
$N_3$  & $32969244672000$     & $59034583740234375$ &$4222124650659840$      
\end{tabular}\vspace*{0.5cm}
\caption{The number of $k$-normal elements of $\mathbb{F}_{q^n}$
over $\mathbb{F}_q$ with $0 \leq k \leq 3$ and $q=2^4$.}
\label{table6}
\end{table}

\section*{Acknowledgements}
Victor G.L.\ Neumann was partially funded by FAPEMIG APQ-03518-18.

\end{document}